\documentclass[12pt,reqno]{amsart}

\usepackage{setspace}
\usepackage{enumitem}
\usepackage{tikz-cd}
\usepackage{tikz}
\usetikzlibrary{positioning} 

\usepackage[linktocpage=true]{hyperref}
\makeatletter
\@namedef{subjclassname@2020}{%
  \textup{2020} Mathematics Subject Classification}
\makeatother
\newtheorem*{theorem*}{Theorem A}
\newtheorem*{theorem**}{Theorem B}

\usepackage[T1]{fontenc}
\usepackage{ae,aecompl}
\usepackage{graphics,cancel,graphicx}
\usepackage{epsfig,psfrag,manfnt}
\usepackage{mathrsfs}
\usepackage{graphicx,mathabx}
\usepackage{bbm,subfigure,rotating,bm,float,mathdots,wasysym,scalerel}

\newlength{\wdth}

\usepackage{stmaryrd}
\usepackage[all,cmtip]{xy}
\usepackage{amssymb,amsmath,amsfonts,amsthm,color}

\usepackage{scalerel,stackengine}
\stackMath
\newcommand\reallywidehat[1]{%
\savestack{\tmpbox}{\stretchto{%
  \scaleto{%
    \scalerel*[\widthof{\ensuremath{#1}}]{\kern-.6pt\bigwedge\kern-.6pt}%
    {\rule[-\textheight/2]{1ex}{\textheight}}%WIDTH-LIMITED BIG WEDGE
  }{\textheight}% 
}{0.5ex}}%
\stackon[1pt]{#1}{\tmpbox}%
}
\newcommand{\sH}{{{\mathscr{H}}^{\infty}}}

\newcommand{\calP}{\mathcal{P}}

\newcommand{\calS}{\mathcal{S}}

\newcommand{\mC}{\mathbb{C}}
\newcommand{\mD}{\mathbb{D}}

\newcommand{\mN}{\mathbb{N}}

\newcommand{\mR}{\mathbb{R}}

\newcommand{\bbz}{\bm{z}}

\newtheorem{theorem}{Theorem}[section]
\newtheorem{lemma}[theorem]{Lemma}
\newtheorem{corollary}[theorem]{Corollary}
\newtheorem{proposition}[theorem]{Proposition}
\newtheorem{conjecture}[theorem]{Conjecture}

\theoremstyle{definition}

\theoremstyle{definition}

\theoremstyle{definition}

\theoremstyle{definition}
\newtheorem{example}[theorem]{Example}

\begin{document}



\title[]{On Banach subalgebras of $\mathscr{H}^\infty$ \\consisting of lacunary Dirichlet series
\vspace{-0.12cm}}

\author[]{Amol Sasane}
\address{Department of Mathematics \\London School of Economics\\
     Houghton Street\\ London WC2A 2AE\\ United Kingdom}
\email{A.J.Sasane@LSE.ac.uk}
 
\maketitle

\vspace{-1cm}
   
\begin{abstract}
\begin{spacing}{1.35}
Let $\mathscr{H}^\infty$ be the set of all Dirichlet series 
  $
\textstyle 
f\!=\!{\scaleobj{0.81}{\sum\limits_{n=1}^\infty}} a_nn^{-s}$ (where $a_n\!\in\! \mC$ for all $n\in \mN=\{1,2,3,\cdots\}$) that converge at each $s$ in $\mC_0\!=\!\{s\in \mC\!:\! \text{Re}(s)\!>\!0\}$, 
such that $\|f\|_{\infty}\!=\!\sup_{s\in \mC_0}|f(s)|\!<\!\infty$. Then $\mathscr{H}^\infty$ is a Banach algebra with pointwise operations and the supremum norm $\|\cdot\|_\infty$, and has been studied in earlier works. The article introduces a new family of Banach subalgebras $\mathscr{H}^\infty_{S}$ of $\mathscr{H}^\infty$. For  $S\!\subset\! \mN$,  
let $\mathscr{H}^\infty_{S}$ be the set of all  elements
    $
\textstyle 
%f\!=\!
{\scaleobj{0.81}{\sum\limits_{n=1}^\infty}} a_nn^{-s}\in \mathscr{H}^\infty
$
such that for all  $n\in \mN\setminus S$,  $a_n\!=\!0$. Then $\mathscr{H}^\infty_{S}$ is a unital Banach subalgebra of $\mathscr{H}^\infty$ with the supremum norm if and only if $S$ is a multiplicative subsemigroup of $\mN$ containing $1$. 
It is  shown that for such $S$, $\mathscr{H}^\infty_{S}$ is the multiplier algebra of $\mathscr{H}^2_S$, where 
$\mathscr{H}^2_S$ is the Hilbert space of all 
%Dirichlet series 
$ 
\textstyle 
f\!=\!{\scaleobj{0.81}{\sum\limits_{n\in S}}} a_nn^{-s}$ such that $\|f\|_2\!:=\!({\scaleobj{0.81}{\sum\limits_{n\in S}}} |a_n|^2)^{\frac{1}{2}}\!<\!\infty$. 
 A characterisation of the  group of units in $\mathscr{H}^\infty_{S}$ is given, by showing an analogue of the Wiener $1/f$ theorem for $\mathscr{H}^\infty_{S}$. 
 If $S$ has a set of generators allowing a unique representation of each element of $S$, then it is shown that the Bass stable rank of $\mathscr{H}^\infty_S$ is infinite.
% The image of $\mathscr{H}^\infty_S$ under the Bohr transform is the subalgebra of the Hardy algebra on the unit ball of $c_0$ with certain vanishing derivatives at $\bm{0}$.
\end{spacing}
\end{abstract}

\vspace{-0.6cm}

\section{Introduction}

\noindent The aim of this article is 
to introduce and study some algebraic-analytic properties of a particular family $\mathscr{H}^\infty_S$ (defined in \S\ref{subsec_subalgebra_S}) of Banach algebras 
that are contained in the Hardy algebra $\mathscr{H}^\infty$ of Dirichlet series (recalled in \S\ref{subsec_Hardy_Dirichlet}). The  motivation is twofold: there has been old and recent interest in studying various 
Banach algebras of Dirichlet series (see, e.g., \cite{BruSas2}, \cite{HedLinSei}), 
and the Banach algebras $\mathscr{H}^\infty_S$ we study are also the `Dirichlet series analogue' of the previously studied (see, e.g., \cite{DPRS}) Banach subalgebra 
 $
\textstyle 
H^\infty_1\!=\!\{f\in H^\infty:f'(0)\!=\!0\}
$ 
   of  $H^\infty$ (the classical Hardy algebra consisting of bounded and holomorphic functions on $\mD:=\{z\in \mC:|z|<1\}$, 
with pointwise operations and the supremum norm, $\|f\|_\infty=\sup_{z\in \mD} |f(z)|$ for $f\in H^\infty$).

\subsection{The Banach algebra $\mathscr{H}^\infty$} 
\label{subsec_Hardy_Dirichlet}
Define $\mC_0\!=\!\{s\in \mC: \text{Re}(s)>0\}$. 
Let $\mathscr{H}^\infty$ be the set of all Dirichlet series 
 $$
\textstyle 
f\!=\!{\scaleobj{0.81}{\sum\limits_{n=1}^\infty}}
a_nn^{-s}
$$  
(where $a_n\in \mC$ for each $n\in \mN$) that converge at each $s\in \mC_0$, 
such that 
$$
\textstyle \|f\|_{\infty}:=\sup\limits_{s\in \mC_0}|f(s)|<\infty.
$$ 
We call $a_n$ the {\em $n^{\text{th}}$ coefficient of $f$}. 
With pointwise operations and the supremum norm, 
$\mathscr{H}^\infty$ is a Banach algebra (introduced in \cite{HedLinSei}). 
Multiplication in $\mathscr{H}^\infty$ is also given by the formula: 
$$
\textstyle 
({\scaleobj{0.81}{\sum\limits_{n=1}^\infty}}
a_n n^{-s}
) 
\cdot ({\scaleobj{0.81}{\sum\limits_{n=1}^\infty}} 
b_nn^{-s}) 
= {\scaleobj{0.81}{\sum\limits_{n=1}^\infty}}
 ({\scaleobj{0.81}{\sum\limits_{(\mN\owns)\;\!d\mid n} }} \;\!a_d b_{\frac{n}{d}})  
n^{-s},
$$
where the notation $d\mid n$ means the integer $d$ divides the integer $n$. 
The unit element is $\mathbf{1}:=\sum_{n=1}^\infty \delta_{n1}n^{-s}$, 
where $\delta_{n1}=0$ for $n\neq 1$ and $ \delta_{11}=1$. 

In
\cite[Theorem~3.1]{HedLinSei}, it was shown that the Banach algebra
$\sH$ is precisely the multiplier space of the Hilbert space $\mathscr{H}^2$, 
where 
$$
\textstyle 
\mathscr{H}^2=\Big\{f\!=\!{\scaleobj{0.81}{\sum\limits_{n=1}^\infty}} 
a_nn^{-s}:\text{such that }
\textstyle  
\|f\|_2\!:=\!\sqrt{{\scaleobj{0.81}{\sum\limits_{n=1}^\infty}} |a_n|^2}\!<\!\infty\Big\}.
$$
The inner product in $\mathscr{H}^2$ is given by 
$$
\textstyle 
\langle f,g\rangle\!=\!{\scaleobj{0.81}{\sum\limits_{n=1}^\infty}} a_n\overline{b_n}, \text{ where }
f\!=\!{\scaleobj{0.81}{\sum\limits_{n=1}^\infty}} 
a_nn^{-s} 
\text{ and }
g\!=\!{\scaleobj{0.81}{\sum\limits_{n=1}^\infty}} 
b_nn^{-s}\in \mathscr{H}^2.
$$
Each element $f\in \mathscr{H}^2$ defines a holomorphic function in the half-plane $\mC_{{\scaleobj{0.81}{\frac{1}{2}}}}\!=\!\{s\in \mC: \text{Re}\;\!s\!>\!\frac{1}{2}\}$. 
If $\zeta$ denotes the Riemann zeta function, 
$$
\textstyle \zeta(s)={\scaleobj{0.81}{\sum\limits_{n=1}^\infty}} 
n^{-s}, \quad 
\text{Re}\;\!s\!>\!1,
$$
then for each $a\in \mC$ such that $\text{Re}\;\!a\!>\!\frac{1}{2}$, the translate 
 $$
 \textstyle 
 \zeta_a(s)={\scaleobj{0.81}{\sum\limits_{n=1}^\infty}} 
{\scaleobj{0.81}{\cfrac{\raisebox{-0.24em}{\text{$1$}}}{\raisebox{0.24em}{\text{$n^a$}}}}}\;\!n^{-s}
$$
belongs to $\mathscr{H}^2$. 
For 
$$
\textstyle f={\scaleobj{0.81}{\sum\limits_{n=1}^\infty}} 
a_nn^{-s}\in \mathscr{H}^\infty,
$$ 
we have 
$\|f\|_2\!\le\! \|f\|_\infty$ (so that $\mathscr{H}^\infty\!\subset\! \mathscr{H}^2$), and in particular, 
\begin{equation}
\label{DGMS_ineq}
\textstyle
|a_n|\le \|f\|_\infty\text{ for all }n\in \mN
\end{equation}
 (see, e.g., \cite[Proposition~1.19]{DGMS}).  
  The set 
 $
\textstyle 
\{e_n:=n^{-s}: n\in \mN\}
$ 
forms an orthonormal basis for $\mathscr{H}^2$. For $a\in \mC_{{\scaleobj{0.81}{\frac{1}{2}}}}$, with
$$
\textstyle 
K_a(s):={\scaleobj{0.81}{\sum\limits_{n=1}^\infty}}  e_n(s) \overline{e_n(a)}=\zeta_{\overline{a}}(s)=\zeta(s+\overline{a}),
$$
 we have $f(a)=\langle f,K_a\rangle$ for all elements $f\in \mathscr{H}^2$. 
The Hilbert space $\mathscr{H}^2$  is a reproducing kernel Hilbert space with  kernel function given by  $K_{\mathscr{H}^2}(s,a)=\zeta(s+\overline{a})$, for $s,a\in \mC_{{\scaleobj{0.81}{\frac{1}{2}}}}$.

\subsection{The set $\mathscr{H}^\infty_S$}
\label{subsec_subalgebra_S} 

  For $S\subset \mN$, define $\mathscr{H}^\infty_{S}$ to be the set of all elements 
$$
\textstyle 
f={\scaleobj{0.81}{\sum\limits_{n=1}^\infty}} a_nn^{-s} \in \mathscr{H}^\infty
$$
such that for all $n\in\mN\setminus  S$, $a_n\!=\!0$. 

%{\small  \tableofcontents }

\subsection{Organisation of the article} 
We show in Section~\ref{section_2} that $\mathscr{H}^\infty_S$ is a unital Banach subalgebra of $\mathscr{H}^\infty$ with the supremum norm if and only if $S$ is a multiplicative subsemigroup of $\mN$ containing $1$. In Section~\ref{section_3}, we show that just as $\mathscr{H}^\infty$ is exactly the multiplier algebra of $\mathscr{H}^2$, the Banach algebra $\mathscr{H}^\infty_S$ is exactly the multiplier algebra of a certain Hilbert subspace $\mathscr{H}^2_S$ of $\mathscr{H}^2$. In Section~\ref{section_4}, we characterise the group of units in $\mathscr{H}^\infty_S$. In Section~\ref{section_5} we relate $\mathscr{H}^\infty_S$ to a natural Banach subalgebra of the Hardy algebra $H^\infty(B_{c_0})$ on the unit ball $B_{c_0}$ of $c_0$ (space of complex sequences converging to $0$ with termwise operations and the supremum norm), with vanishing derivatives of certain orders at $\bm{0}\in c_0$. Finally, in Section~\ref{section_6}, we prove that if $S$ has a set of generators allowing a unique representation of each element of $S$, then the Bass stable rank of $\mathscr{H}^\infty_S$ is infinite, and also state a related conjecture.

\section{When is $\mathscr{H}^\infty_S$ an algebra?}
\label{section_2}

\noindent We show that $\mathscr{H}^\infty_{S}$ is a unital Banach subalgebra of $\mathscr{H}^\infty$ with the supremum norm if and only if $S$ is a multiplicative subsemigroup of $\mN$ containing $1$. Recall that a subset $\calS$ of a semigroup $\Sigma$ is a {\em subsemigroup} of $\Sigma$ if $\calS\calS:=\{s_1s_2:s_1,s_2\in \calS\}\subset \calS$. By a subalgebra $B$ of a complex algebra $A$ (see, e.g., \cite[Definition~10.1]{Rud}), we simply mean that $B$ is closed under the algebraic operations (addition, scalar multiplication and multiplication) inherited from $A$. 

\begin{proposition}
\label{prop_1.1}
The following are equivalent$\;\!:$

\noindent {\em (1)}  $S$ is a multiplicative subsemigroup of $\mN$. 

\noindent {\em (2)} $\mathscr{H}^\infty_{S}$ is a subalgebra of $\mathscr{H}^\infty$. 

\noindent {\em (3)}  $\mathscr{H}^\infty_{S}$ is a Banach algebra with the supremum norm. 
\end{proposition}
\begin{proof}(1)$\Rightarrow$(2): 
 Suppose $S$ is a multiplicative subsemigroup of $\mN$. It is clear that $\mathscr{H}^\infty_{S}$ is nonempty since the zero function $\bm{0}\in \mathscr{H}^\infty$ belongs to $\mathscr{H}^\infty_{S}$. 
Closure under addition and scalar multiplication is obvious. We show closure under multiplication. 
Let 
$$
\textstyle 
f={\scaleobj{0.81}{\sum\limits_{n=1}^\infty }} a_nn^{-s}\in \mathscr{H}^\infty_S \text{ and }g={\scaleobj{0.81}{\sum\limits_{n=1}^\infty}} b_nn^{-s} \in \mathscr{H}^\infty_S.
$$
  The $n^{\text{th}}$ coefficient of $fg$ is given by 
 \begin{equation}
 \label{sum_for c_n}
 \textstyle 
  c_n:={\scaleobj{0.81}{\sum\limits_{d\mid n} }} \;\!a_d b_{\frac{n}{d}}.
  \end{equation}
  Let $n\not\in S$. We claim $c_n\!=\!0$. 
  We have the following two cases. 
  
  \vspace{0.03cm}
  
  \noindent $1^\circ$ Let $d\not\in S$. So $a_d\!=\!0$. In the 
  sum \eqref{sum_for c_n} for $c_n$, the summand $a_d b_{\frac{n}{d}}\!=\!0$. 
  
  \vspace{0.03cm}
  
  \noindent $2^\circ$ Let $d\in S$. If  $\frac{n}{d}\in  S$ too,   then $n=d\cdot \frac{n}{d}\in S$ (as $S$ is a multiplicative \noindent \phantom{$1^\circ$ }semigroup), a contradiction to $n\not\in S$. Thus $\frac{n}{d}\not\in S$. Hence  $b_{\frac{n}{d}}=0$, 
  \noindent \phantom{$1^\circ$ }and so in the sum \eqref{sum_for c_n} for $c_n$, the summand $a_d b_{\frac{n}{d}}=0$ again. 
  
  \vspace{0.03cm}
 
 \noindent 
Thus $c_n\!=\!0$. So $fg\in \mathscr{H}^\infty_{S}$. Consequently, $\mathscr{H}^\infty_{S}$ is a subalgebra of $\mathscr{H}^\infty$.

\vspace{0.12cm}

\noindent (2)$\Rightarrow$(1): Suppose $\mathscr{H}^\infty_S$ is a subalgebra of $\mathscr{H}^\infty$. Let   $m,n\in S$. 
Define $\textstyle 
f\!=\!n^{-s}$ and  $g\!=\!m^{-s}$. 
Then $f,g\in \mathscr{H}^\infty_S$. As $\mathscr{H}^\infty_S$ is an algebra with pointwise operations,  $fg\in \mathscr{H}^\infty_S$. 
But $fg\!=\!(nm)^{-s}$. 
%So it cannot be the case that $nm\not\in S$.  
Hence $nm\in S$.  Consequently, $S$ is a multiplicative subsemigroup of $\mN$.

\vspace{0.12cm}

\noindent As (3)$\Rightarrow$(2) is trivial, it remains to show  
 (2)$\Rightarrow$(3), i.e., $\mathscr{H}^\infty_S$ is a closed subset of $\mathscr{H}^\infty$. Let $(f_m)_{m\in \mN}$ be a Cauchy sequence in $\mathscr{H}^\infty_S$. Write 
 $$
  \textstyle
 f_m={\scaleobj{0.81}{\sum\limits_{n=1}^\infty}} \;\!a^{(m)}_n n^{-s}.
 $$
 Then $(f_m)_{m\in \mN}$ is also  Cauchy in the Banach algebra $\mathscr{H}^\infty$, and hence it converges to some element
 $$
 \textstyle
 f={\scaleobj{0.81}{\sum\limits_{n=1}^\infty}}\;\!
 a_n n^{-s} \in \mathscr{H}^\infty.
 $$
 Let $n\in \mN\setminus S$. Then $a_n^{(m)}\!=\!0$ for all $m\in \mN$. Using \eqref{DGMS_ineq}, we get
  $$
 \textstyle |a_n|=|a_n-0|=|a_n-a_n^{(m)}|\le \|f-f_m\|_\infty\;\text{ for all }m\in \mN.
 $$ 
  Passing to the limit as $m\to \infty$,  $|a_n|\!\le\! 0$, i.e., $a_n\!=\!0$. So $f\in \mathscr{H}^\infty_S$. 
 \end{proof}
 
 \begin{proposition}
 Let $S$ be a multiplicative subsemigroup of $\mN$. 
 
 \noindent Then $\mathscr{H}^\infty_S$ is unital if and only if $1\in S$. 
 \end{proposition}
 \begin{proof} If $1\in S$, then $\bm{1}\in \mathscr{H}^\infty_S$, and so $\mathscr{H}^\infty_S$ is unital. 
  
  \noindent Conversely, if $\bm{1}\in \mathscr{H}^\infty_S$, then it cannot be the case that $1\not\in S$. 
 \end{proof}

\begin{example}
\label{examples}$\;$

\noindent 
Some examples of multiplicative subsemigroups of $\mN$ containing $1$ are listed below.
 
 \smallskip 
 
\noindent (1)  For a prime $p$, define $S\!=\!\{p^k: k\in \mN\cup\{0\}\}$.  

\noindent  \phantom{(1) }(Then  
 
 \vspace{-0.12cm}
 
 $\phantom{AAAAAAAAi}
\textstyle \mathscr{H}^\infty_{S}\owns  {\scaleobj{0.81}{\sum\limits_{n=1}^\infty}}
 a_n n^{-s} \mapsto {\scaleobj{0.81}{\sum\limits_{n=1}^\infty }}a_n z^n\in H^\infty
 $
 
  \vspace{0.06cm}
  
 \noindent \phantom{(1) }is a Banach algebra isomorphism from $\mathscr{H}^\infty_{S}$ to $H^\infty$.)
 
 \noindent 
 \phantom{(1) }More generally, for any $n\in \mN$, define $S:=\{n^k: k\in \mN\cup\{0\}\}$. 
 
  \smallskip

\noindent (2) For primes $p_{i_1}\!<\!\cdots\!<p_{i_m}$, let $S\!=\!\{p_{i_1}^{k_1}\cdots p_{i_m}^{k_m}: k_1,\cdots, k_m\in \mN\cup\{0\}\}$. 

\smallskip 

\noindent (3) Let $\chi_0$ be the principal Dirichlet character of modulus $m\in \mN$, 
 $$
 \textstyle 
 \chi_0 (n)\!=\!\Big\{\!\begin{array}{l}
 0\;\text{ if } \;\text{gcd}(n,m)\!>\!1,\\
 1\;\text{ if } \;\text{gcd}(n,m)\!=\!1.
 \end{array}
 $$
 \phantom{(1) }where $\text{gcd}(n,m)$ denotes the greatest common divisor of $n,m\in \mN$. 
 
 \noindent \phantom{(1) }The set  
  $
 \textstyle 
 S_m\!=\!\{n\in \mN: \chi_0(n)\!\neq \!0\}\!=\!\{n\in \mN:\text{gcd}(n,m)\!=\!1\}.
 $  is a 
 \noindent \phantom{(1) }multiplicative subsemigroup of $\mN$ containing $1$.
 Let 
 $p_1\!<\!p_2\!<\!p_3\!<\!\cdots$ 
 \noindent \phantom{(1) }be all the prime numbers arranged in ascending 
order.  For $n\;\!\in\;\! \mN$, 

\noindent \phantom{(1) }the fundamental theorem of arithmetic gives 
 the existence of a 
 \noindent \phantom{(1) }unique compactly supported sequence $(\nu_k(n))_{k\in \mN}$ 
 of  nonnegative in-
\noindent \phantom{(1) }tegers such that 
$$
\textstyle 
n={\scaleobj{0.81}{\prod\limits_{k=1}^\infty}} \; p_k^{\nu_k(n)}.
$$
 \phantom{(1) }As
  $$
 \textstyle 
 \text{gcd}(n,m)\!=\!
 {\scaleobj{0.81}{\prod\limits_{k=1}^\infty}} \; p_k^{\min \{\nu_k(n),\nu_k(m)\}}\!,
 $$
 \noindent \phantom{(1) }we get $\text{gcd}(n,m)\!=\!1$ if and $\!$only if 
 $$
 \textstyle \text{for all }k\!\in\! \mN, \;\min \{\nu_k(n),\nu_k(m)\}\!=\!0.
 $$
  \phantom{(1) }Hence 
  $
 \textstyle
 S_m=\big\{n={\scaleobj{0.81}{\prod\limits_{k=1}^\infty}} \;\!p_k^{\nu_k(n)}: \text{for all }k\in \mN,\; \nu_k(n)=0 \text{ if } \nu_k(m)\!>\!0\big\}.
 $
 
 \noindent 
    \phantom{(1) }As the intersection of multiplicative subsemigroups of $\mN$ is again 
 \phantom{(1) }a multiplicative subsemigroup of $\mN$, for any subset $F\subset \mN$, 
 $$
 \textstyle 
 S_F:={\raisebox{0.15em}{\scaleobj{0.81}{\bigcap\limits_{m\in F}}}} \;\! S_m
 $$
  \phantom{(1) }is a multiplicative subsemigroup of $\mN$ containing $1$. If $F$ is a finite 
  \phantom{(1) }nonempty set  $F=\{m_1,\cdots, m_n\}\subset \mN$, 
  then 
  $$
  \textstyle 
 S_F:={\raisebox{0.15em}{\scaleobj{0.81}{\bigcap\limits_{m\in F}}}} \;\! S_m=S_{\text{lcm}(m_1,\cdots, m_n)},
 $$
 \phantom{(1) }where $\text{lcm}(m_1,\cdots, m_n)$ is the least common multiple of $m_1,\cdots, m_n$. 
 
 \smallskip 
 
 \noindent (5) $S=\mN$. (Then $\mathscr{H}^\infty_S=\mathscr{H}^\infty$.) 
 
 \smallskip 
 
 \noindent (6) $S=\{1\}$. (Then $\mathscr{H}^\infty_S$ is isomorphic to the  Banach algebra $\mC$.)
 
 \smallskip 
 
 \noindent (7) Given $m\in \mN$, $S:=\{n^m:n\in \mN\}$. 
 
 \noindent \phantom{(1) }(For example, if $m\!=\!2$, then $S$ is 
the set of all 
perfect squares.)
 
 \smallskip 
 
 \noindent (8) Let $S=\{1\}\cup\{n\in \mN: \text{there exist } x,y\in \mN\text{ such that }n\!=\!x^2+y^2\}$. 
 
 \noindent \phantom{(1) }Then $S$ is a  multiplicative subsemigroup of $\mN$ containing $1$. 
 
 \noindent \phantom{(1) }The sum of two squares theorem (see, e.g., \cite[Thm.~7, \S3, Chap.~IV]{Cha}) 
 \noindent \phantom{(1) }gives an alternative description of $S$: $n\in S$ if and only if  all prime 
 \noindent \phantom{(1) }factors of $n$ of the form $4k + 3$ ($k\!\in\! \mN\!\cup\! \{0\}$) appear with an even \noindent \phantom{(1) }exponent in the prime factorisation of $n\in \mN$.
 \hfill$\Diamond$
 \end{example}
 
 \section{$\mathscr{H}^\infty_S$ as the multiplier algebra of $\mathscr{H}^2_S$}
 \label{section_3}
 
 \noindent It was shown in \cite{HedLinSei} that $\mathscr{H}^\infty$ is exactly the multiplier algebra of $\mathscr{H}^2$, 
 that is, a function $f$ defined on $\mC_{{\scaleobj{0.81}{\frac{1}{2}}}}$ satisfies $fg\in \mathscr{H}^2$ for all $g\in \mathscr{H}^2$ if and only if $f$ has an extension to $\mC_0$ which is an element of $\mathscr{H}^\infty$. 
 Moreover, 
  $
 \textstyle 
 \|f\|_{\infty}=\sup\limits_{g\in \mathscr{H}^2,\;\! \|g\|_2\le 1}\|fg\|_2.
 $
 
 Analogously, in this section, we will show that more generally, $\mathscr{H}^\infty_S$ is exactly the multiplier algebra of $\mathscr{H}^2_S$. We first define $\mathscr{H}^2_S$. 
 
 For any subset $S\subset \mN$, we define 
 $
 \mathscr{H}^2_S$ to be the set of all 
 $$
 \textstyle 
 f={\scaleobj{0.81}{\sum\limits_{n=1}^\infty}} a_n n^{-s}  \in \mathscr{H}^2
 $$
 such that for all $n\in \mN\setminus S$,  $a_n\!=\!0$. Then $ \mathscr{H}^2_S$ is a closed subspace of $\mathscr{H}^2$, and 
  $
 \textstyle
 \{n^{-s}: n\in S\}
 $ 
 forms an orthonormal basis for $\mathscr{H}^2_S$. Define lacunary zeta function $\zeta_S$ by 
 $$
 \textstyle
 \zeta_S(s)=
 {\scaleobj{0.81}{\sum\limits_{n\in  S}}} \;\!n^{-s} , 
 \quad 
\text{Re}\;\!s\!>\!1.
$$
 For $a\in \mC_{{\scaleobj{0.81}{\frac{1}{2}}}}$,
$$
\textstyle 
{\scaleobj{0.81}{\sum\limits_{n\in S}}}  e_n(s) \overline{e_n(a)}=\zeta_S(s+\overline{a}),
$$
 and we have $f(a)=\langle f,\zeta_S(\cdot+\overline{a})\rangle$ for all $f\in \mathscr{H}^2_S$. 
The Hilbert space $\mathscr{H}^2_S$  is a reproducing kernel Hilbert space with  kernel function given by  $K_{\mathscr{H}^2_S}(s,a)=\zeta_S(s+\overline{a})$ for $s,a\in \mC_{{\scaleobj{0.81}{\frac{1}{2}}}}$. 
In particular, in Example~\ref{examples}(2), $\mathscr{H}^\infty_{S_m}$ is a reproducing kernel Hilbert space with the kernel given by 
$K_{\mathscr{H}^\infty_{S_m}}(s,a)=L(s+\overline{a},\chi_0)$, where $L$ is the following Dirichlet $L$-series:
$$
\textstyle 
L(s,\chi_0)={\scaleobj{0.81}{\sum\limits_{n=1}^\infty}} \chi_0(n)n^{-s},\quad \text{Re}\;\!s\!>\!1.
$$
(Peripherally, a natural question is:  
Is there a characterisation of the multiplicative subsemigroups $S$ of $\mN$ containing $1$, for which the  lacunary zeta functions $\zeta_S$ arise from modular forms? See, e.g., \cite{Apo}, for background on modular forms and their link to Dirichlet series.)

 We have the following:
 
\begin{proposition}$\;$
\label{prop_14_8_2025_1801}

\noindent 
Let $S$ be a multiplicative subsemigroup of $\mN$ containing $1$. 
Then $\mathscr{H}^\infty_S$ is exactly the multiplier algebra of $\mathscr{H}^2_S,$ that is$,$ a function $f$ defined on $\mC_{{\scaleobj{0.81}{\frac{1}{2}}}}$ satisfies $fg\in \mathscr{H}^2_S$ for all $g\in \mathscr{H}^2_S$ if and only if $f$ has an extension to $\mC_0$ which is an element of $\mathscr{H}^\infty_S$. Moreover$,$ 
  $
 \textstyle 
 \|f\|_{\infty}=\sup\limits_{g\in \mathscr{H}^2_S,\;\! \|g\|_2\le 1}\|fg\|_2.
 $ 
\end{proposition}

\noindent We will show this along the same lines as the proof for $\mathscr{H}^\infty$-$\mathscr{H}^2$ case given in \cite[Theorem~6.4.7]{QueQue}. 

\noindent Let $\calP$ denote the set of all {\em Dirichlet polynomials}, that is,
$$
\textstyle
\calP=\Big\{f\!=\!{\scaleobj{0.81}{\sum\limits_{n=1}^\infty}} a_n n^{-s}: 
{\scaleobj{0.9}{\begin{array}{ll}\text{there exists an }N\in \mN\text{ such that}\\
\text{for all }n\!>\!N, \; a_n\!=\!0\end{array}}}\Big\}.
$$
For $r\in \mN$, let $N_r\!=\!\{n\!=\!p_1^{k_1}\cdots p_r^{k_r}: k_1,\cdots, k_r \in \mN\cup\{0\}\}$. 
Define
$$
\textstyle
\calP_r=\big\{f={\scaleobj{0.81}{\sum\limits_{n\in N_r}}} a_n n^{-s}\in \calP: \text{for all }n\in N_r, \;\! a_n\in \mC\big\}.
$$
We recall  \cite[Lemma~6.4.9]{QueQue}, and include its short proof (correcting a  typo: 
the definition of $C_{s,r}$ should include a missing square root).   

\begin{lemma}
\label{lemma_14_8_2025_1840}
For all $s\in \mC_0,$ and for all $r\in \mN,$ there exists a constant $C_{s,r}\!>\!0$ such that for all $f\in \calP_r,$ we have $|f(s)|\le C_{s,r}\|f\|_2$. 
\end{lemma}
\begin{proof} Let $\sigma\!>\!0$, $t\in \mR$ and $s\!:=\!\sigma+it \in \mC_0$. Let 
$$
\textstyle f\!=\!{\scaleobj{0.81}{\sum\limits_{n\in N_r}}} a_n n^{-s}\in \calP_r.
$$
  Define $C_{s,r}\!>\!0$ by 
 $$
\textstyle 
C_{s,r}^2:= {\scaleobj{0.81}{\sum\limits_{n\in N_r}}} n^{-2\sigma}
={\scaleobj{0.81}{\prod\limits_{j=1}^{r}}}(1-p_j^{-2\sigma})^{-1} <\infty,
$$ 
where the second equality follows from the  truncated Euler identity (obtained by expanding each factor as the sum of a geometric series).
The Cauchy-Schwarz inequality gives 
$$
\textstyle 
\phantom{AAAAAAA}
|f(s)|\le ({\scaleobj{0.81}{\sum\limits_{n\in N_r}}} |a_n|^2)^{\frac{1}{2}} ({\scaleobj{0.81}{\sum\limits_{n\in N_r}}} n^{-2\sigma})^{\frac{1}{2}} =C_{s,r}\|f\|_2.
\phantom{AAAAAAA}
\qedhere
$$
\end{proof}

\noindent For $\varphi \in \text{L}^1(\mR)$, let $\widehat{\varphi}$ denote the Fourier transform of $\varphi$:
$$
\textstyle 
\widehat{\varphi}(\xi)=\int_{\mR} \varphi(t)e^{-i\xi t} dt, \quad \xi \in \mR.
$$
Let
 $
\textstyle 
E=\{\varphi\in \text{L}^1(\mR): \widehat{\varphi} \text{ has compact support}\}.
$ 
Then given 
$$
\textstyle f={\scaleobj{0.81}{\sum\limits_{n=1}^\infty }} a_nn^{-s} \in \mathscr{H}^2,
$$
 and $\varphi \in E$, 
the following `vertical convolution identity' holds:
\begin{equation}
\label{Saksman}
\textstyle 
 {\scaleobj{0.81}{\sum\limits_{n=1}^\infty }} a_n \widehat{\varphi}(\log n) n^{-s}
 =
 \int_\mR f(s+it)\varphi(t) dt , \quad s\in \mC_{{\scaleobj{0.81}{\frac{1}{2}}}}.
 \end{equation}
 (See, e.g.,  \cite[Proof of Theorem 6.4.7]{QueQue}.)

\begin{proof}[Proof of Proposition~\ref{prop_14_8_2025_1801}:] 
If $f\in \mathscr{H}^\infty_S$, and $g\in \mathscr{H}^2_S$, then $ fg \in \mathscr{H}^2$, 
and $\|\varphi f\|_2\le \|\varphi\|_\infty \|f\|_2$.  Let 
$$
\textstyle 
f\!=\!{\scaleobj{0.81}{\sum\limits_{n=1}^\infty}}  a_nn^{-s} 
\in \mathscr{H}^\infty_S\; \text{ and }\;g\!=\!{\scaleobj{0.81}{\sum\limits_{n=1}^\infty}}
b_nn^{-s} \in \mathscr{H}^2_S.
$$
  The $n^{\text{th}}$ coefficient of $fg$ is given by 
 \begin{equation}
 \label{sum_for c_n}
 \textstyle 
  c_n\!:=\!{\scaleobj{0.81}{\sum\limits_{d\mid n} }} \;\!a_d b_{\frac{n}{d}}.
  \end{equation}
  Let $n\not\in S$. We claim that $c_n\!=\!0$. 
  We have the following two cases. 
  
  \vspace{0.03cm}
  
  \noindent $1^\circ$ Let $d\not\in S$. So $a_d\!=\!0$. In the 
  sum \eqref{sum_for c_n} for $c_n$, the summand $a_d b_{\frac{n}{d}}\!=\!0$. 
  
  \vspace{0.03cm}
  
  \noindent $2^\circ$ Let $d\in S$. If  $\frac{n}{d}\in  S$ too,   then $n=d\cdot \frac{n}{d}\in S$ (as $S$ is a multiplicative \noindent \phantom{$1^\circ$ }semigroup), a contradiction to $n\not\in S$. Thus $\frac{n}{d}\not\in S$. Hence  $b_{\frac{n}{d}}=0$, 
  \noindent \phantom{$1^\circ$ }and so in the sum \eqref{sum_for c_n} for $c_n$, the summand $a_d b_{\frac{n}{d}}=0$ again. 
  
  \vspace{0.03cm}
 
 \noindent 
Thus $c_n\!=\!0$. So $fg\in \mathscr{H}^2_{S}$. Hence if $f\in \mathscr{H}^\infty_S$, then the multiplication map  $\mathscr{H}^2_S\owns g\mapsto M_f g:= fg \in \mathscr{H}^2_S$ is well-defined, and $\|M_f\|\le \|f\|_\infty$. 

\smallskip

Next, suppose  $f:\mC_{{\scaleobj{0.81}{\frac{1}{2}}}}\to \mC$ is such that $fg\in \mathscr{H}^2_S$ for all $g\in\mathscr{H}^2_S$. 
Let $M_f:\mathscr{H}^2_S\to \mathscr{H}^2_S$ be the linear map of pointwise multiplication by $f$.  
As $\bm{1}\in \mathscr{H}^2_S$, we have $f\!=\!M_f(\bm{1})\in \mathscr{H}^2_S$. By the closed graph theorem, $M_f$ is a bounded operator. Denote the operator norm of $M_f$ by $\|M_f\|$. Let 
$$
\textstyle 
f={\scaleobj{0.81}{\sum\limits_{n=1}^\infty}} a_n n^{-s}, \quad s\in \mC_{{\scaleobj{0.81}{\frac{1}{2}}}}.
$$

\noindent {\bf Step 1.} Suppose first that $f$ is a Dirichlet polynomial. We claim that $\|f\|_\infty=\|M_f\|$. 
Fix $r\in \mN$ such that $f\in \calP_r$ and let $s\in \mC_0$. 
We argue by induction that for all $k\in \mN$, we have $\|f^k\|_2\le \|M_f\|_k$. If $k=1$, then 
$\|f\|_2^1\!=\!\|M_f(\bm{1})\|_2\!\le\! \|M_f\|\|\bm{1}\|_2\!=\!\|M_f\| 1\!=\!\|M_f\|^1$. 
If $\|f^k\|_2\!\le \!\|M_f\|^k$ for a $k\in \mN$, then  
$\|f^{k+1}\|_2\!=\!\|M_f(f^k)\|_2\!\le\! \|M_f\|\|f^k\|_2\!\le\! \|M_f\|\|M_f\|^k\!=\!\|M_f\|^{k+1}$, completing the induction step. Lemma~\ref{lemma_14_8_2025_1840} applied to $f^k\in \calP_r$ implies for $s\in \mC_0$ that 
$$
\textstyle 
|f(s)|^k \le C_{s,r} \|f^k\|_2\le C_{s,r}\|M_f\|^k,
$$
and so $|f(s)|\le C_{s,r}^{\frac{1}{k}} \|M_f\|$. Passing to the limit that $k\to \infty$ now yields $|f(s)|\!\le\! 1\|M_f\|$. 
As $s\in \mC_0$ was arbitrary, $\|f\|_\infty\! \le\!  \|M_f\|$. Also, as $f\in \mathscr{H}^2_S$ is Dirichlet polynomial and $\|f\|_\infty\! <\! \infty$, we have that $f\in \mathscr{H}^\infty_S$. But then $f$ is a multiplier on $\mathscr{H}^2_S$ and $\|M_f\|\! \le\!  \|f\|_\infty$ by the first part of the proof. Consequently, $\|f\|_\infty\! =\! \|M_f\|$, as  claimed. 

\smallskip 

\noindent {\bf Step 2.} Now consider the general case when $f$ need not be Dirichlet polynomial. 
For a function $\varphi\in E$, we define 
$$
\textstyle  
P_{\varphi}(s)={\scaleobj{0.81}{\sum\limits_{n=1}^\infty}} a_n \widehat{\varphi}(\log n) n^{-s}.
$$
As $\widehat{\varphi}$ has compact support, and $\log n\to \infty$ as $n\to \infty$, it follows that $P_\varphi$ is a Dirichlet polynomial. We claim that $\|M_{P_\varphi}\|\le \|M_f\|\|\varphi\|_1$. For a $t\in \mR$, define   the vertical translation operator $T_t$ by  $(T_t g)(s)=g(s+it)$ for all $g\in \mathscr{H}^2_S$. Then $T_t:
\mathscr{H}^2_S\to \mathscr{H}^2_S$ is a linear isometry on $\mathscr{H}^2_S$, and $T_t f$ is a multiplier on $\mathscr{H}^2_S$ satisfying $\|M_{T_t f}\|=\|M_f\|$. Indeed, for all $g\in \mathscr{H}^2_S$, we have $(T_t f) g \!=\!T_t (f (T_{-t} g))$, and 
 $$
\textstyle 
\| (T_t f)g\|_2\!=\!\|T_t (f (T_{-t} g))\|_2=\|f(T_{-t} g)\|_2\!\le\! \|M_f\|\|T_{-t} g\|_2\!=\!\|M_f\|\|g\|_2,
$$
giving $\|M_{T_t f}\|\! \le\!  \|M_f\|$. Then also  $\|M_f \|\! =\! \|M_{T_{-t}(T_t f)}\|\! \le \! \|M_{T_t f}\|$. 
The vertical convolution formula \eqref{Saksman}
 yields  for $s\in \mC_{{\scaleobj{0.81}{\frac{1}{2}}}}$ and $g\in \mathscr{H}^2_S$ that:
 $$
 \textstyle 
 \begin{array}{rcl}
 (P_{\varphi}g)(s)\!\!\!\!&=&\!\!\!(\int_\mR f(s +it)\varphi (t) dt) g(s)= \int_\mR (T_t f)(s) g(s) \varphi(t) dt\\[0.21cm]
 \!\!\!\!&=&\!\!\!\int_\mR ((T_t f)g)(s) \varphi(t) dt .
 \end{array}
 $$
 We have $P_{\varphi}g=\int_\mR ((T_t f)g)(\cdot) \varphi(t) dt$ in $\mathscr{H}^2_S$, where the right-hand side is a vector-valued Pettis integral in $\mathscr{H}^2_S$, and 
 $$
 \textstyle 
 \begin{array}{rcl}
 \|M_{P_{\varphi}}g\|_2=\|P_{\varphi}g\|_2\!\!\!\!&\le&\!\!\! \int_\mR \|(T_t f)g\|_2|\varphi(t) |dt \\[0.21cm]
 \!\!\!\!&\le&\!\!\! \int_\mR \|M_{T_t f}\|\|g\|_2|\varphi(t) |dt =\int_\mR \|M_f\|\|g\|_2 |\varphi(t) |dt \\[0.21cm]
 \!\!\!\!&=&\!\!\!\|M_f\|\|g\|_2\int_\mR |\varphi(t) |dt =\|M_f\|\|g\|_2\|\varphi\|_1.
 \end{array}
 $$
 Hence $\|M_{P_{\varphi}}\|\le \|M_f\|\|\varphi\|_1$.
 
 \smallskip 

\noindent {\bf Step 3.} Define the sequence $(\varphi_m)_{m\in \mN}$ in $\text{L}^1(\mR)$  by
$$
\textstyle 
\varphi_m(t)={\scaleobj{0.81}{\cfrac{m}{2\pi}}}({\scaleobj{0.81}{\cfrac{\sin \frac{mt}{2}}{\frac{mt}{2}}}})^2, 
\quad t\in \mR.
$$
Then 

\vspace{-0.3cm}

$\phantom{AAAAAAAAAAAA}
\textstyle 
\widehat{\varphi_m}(\xi)=\max\{1\!-\!{\scaleobj{0.81}{\cfrac{\raisebox{-0.21em}{\text{$|\xi|$}}}{\raisebox{0.3em}{\text{$m$}}}}}, 0\}, 
$

\vspace{0.06cm}

\noindent 
and as $\widehat{\varphi_m}$ has compact support, we have that $\varphi_m\in E$ for all $m\in \mN$. 
Moreover, $\varphi_m\!\ge \!0$, and so 
$$
\textstyle
1=\widehat{\varphi_m}(0)=\int_{\mR} \varphi_m(t) e^{i0 t} dt=\int_{\mR} \varphi_m(t) dt=\int_{\mR} |\varphi_m(t)| dt=\|\varphi_m\|_1.
$$
Then 
$$
\textstyle 
P_{\varphi_m}(s)={\scaleobj{0.81}{\sum\limits_{n=1}^\infty}} a_n \widehat{\varphi_m}(\log n) n^{-s}
=\int_\mR (T_t f)(s) \varphi_m(t) dt, \quad s\in \mC.
$$
Then Steps 1 and 2 imply that for all $m\in \mN$, 
$$
\textstyle 
\|P_{\varphi_m}\|_\infty 
=\|M_{P_{\varphi_m}} \|\le \|M_f\|\|\varphi_m\|_1=\|M_f\|1=\|M_f\|.
$$
Taking a subsequence if necessary, one may assume, thanks to Montel's theorem, that $P_{\varphi_m}$ tends to some $F$ uniformly on compact subsets of $\mC_0$, with 
$$
\textstyle 
\|F\|_\infty:=\sup\limits_{s\in \mC_0}|F(s)|\le \|M_f\|.
$$
 Let $\sigma\!>\!\frac{1}{2}$, $t\in\! \mR$, and $s\!=\!\sigma+it \in\! \mC_{{\scaleobj{0.81}{\frac{1}{2}}}}$. By the Cauchy-Schwarz inequality,
$$
\textstyle 
{\scaleobj{0.81}{\sum\limits_{n=1}^\infty}} |a_n|n^{-\sigma} 
\le 
( {\scaleobj{0.81}{\sum\limits_{n=1}^\infty}} |a_n|^2)^{\frac{1}{2}}
( {\scaleobj{0.81}{\sum\limits_{n=1}^\infty}} n^{-2\sigma})^{\frac{1}{2}}<\infty.
$$
Also $\widehat{\varphi_m}(\log n)\!\to\! 1$ as $m\!\to\! \infty$, and $0\!\le\! \widehat{\varphi_m}(\log n)\!\le\! 1$. 
Given $\epsilon\!>\!0$, let $N\in \!\mN$ be such that 
$$
\textstyle 
{\scaleobj{0.81}{\sum\limits_{n=N+1}^\infty}} |a_n|n^{-\sigma}<\frac{\epsilon}{4}.
$$
Let $m_n$, $n\in \{1,\cdots, N\}$ be such that 
$$
\textstyle 
|\widehat{\varphi_{m_n}}(\log n)-1|\le {\scaleobj{0.81}{\cfrac{\raisebox{-0.21em}{\text{$\epsilon$}}}{2N ({\scaleobj{0.96}{\sum\limits_{n=1}^N}} |a_n|n^{-\sigma}+1)}}}.
$$
Then for $m>\max\{m_1,\cdots, m_N\}$, we have 
$$
\textstyle 
\begin{array}{rcl}
|P_{\varphi_m}(s)-{\scaleobj{0.81}{\sum\limits_{n=1}^\infty}} a_n n^{-s}|
\!\!\!\!&=&\!\!\!\!
|{\scaleobj{0.81}{\sum\limits_{n=1}^\infty}} a_n (\widehat{\varphi_{m_n}}(\log n)-1)n^{-s}|
\\[0.21cm]
\!\!\!\!&\le&\!\!\!\!
{\scaleobj{0.81}{\sum\limits_{n=1}^N}}|a_n||\widehat{\varphi_{m_n}}(\log n)-1| n^{-\sigma}
+
{\scaleobj{0.81}{\sum\limits_{n=N+1}^\infty}}|a_n|2 n^{-\sigma}
\\[0.33cm]
\!\!\!\!&\le&\!\!\!\!
N\frac{\epsilon}{2N}\cdot 1+2\frac{\epsilon}{4}=\epsilon.
\end{array}
$$
Thus for each $s\in \mC_{{\scaleobj{0.81}{\frac{1}{2}}}}$, we have 
$$
\textstyle
P_{\varphi_m}(s)\to {\scaleobj{0.81}{\sum\limits_{n=1}^\infty}} a_n n^{-s}\!=\!f(s)\;\text{ as }m\!\to\! \infty.
$$
Hence $f\!=\!F$ on $\mC_{{\scaleobj{0.81}{\frac{1}{2}}}}$. But $f\in \mathscr{H}^2_S$, and so it is a Dirichlet series. 
We have shown that $f$ has a Dirichlet series which converges in $\mC_{{\scaleobj{0.81}{\frac{1}{2}}}}$, and this $f$ admits a bounded holomorphic extension $F$ to $\mC_0$. Thus it follows that $f\in \mathscr{H}^\infty$. As $f\in \mathscr{H}^2_S \cap \mathscr{H}^\infty$, we get $f\in \mathscr{H}^\infty_S$. 
Moreover, $\|f\|_\infty\!=\!\|F\|_\infty\!\le\! \|M_f\|$. Since also $\|M_f\|\!\le\! \|f\|_\infty$, we get $\|f\|_\infty\!=\!\|M_f\|$. 
\end{proof}
  
 \section{Characterisation of the group of units}
  \label{section_4}
  
\noindent In this section we will show that $f\in \mathscr{H}^\infty_S$ is invertible in $\mathscr{H}^\infty_S$ if and only if 
 $
\textstyle
\inf\limits_{s\in \mC_0} |f(s)|>0.
$

Below, for a unital commutative complex Banach $A$, 
we denote by $A^{-1}$ the multiplicative group of all invertible elements of $A$. 
For $\sigma\in \mR$, let $\mC_\sigma\!:=\!\{s\in\mC: \text{Re}\;\!s\!>\!\sigma\}$. 
Recall that for the Dirichlet series  
$$
\textstyle 
D\!=\!{\scaleobj{0.81}{\sum\limits_{n=1}^\infty}} a_n n^{-s},
$$ 
the {\em abscissa of convergence} is defined by 
$$
\textstyle 
\sigma_c(D)\!=\!
\inf\{\sigma \in \mR: D\textrm{ converges in }\mC_\sigma\}
\in [-\infty,\infty].
$$
 Similarly, the 
abscissa of {\em absolute}  convergence of $D$ 
is given by
$$
\textstyle 
\sigma_a(D)
\!=\!
\inf\{\sigma \in \mR: D\textrm{ converges absolutely in } 
\mC_\sigma\}.
$$
We have $-\infty\!\le\!\sigma_c(D)\!\le\!  \sigma_a(D)\!\le\!\infty$. 
Also, $\sigma_a(D)\!\le\! \sigma_c(D)+1$, 
see, e.g., \cite[Proposition~1.3]{DGMS}. 
 
\begin{theorem}
\label{te2.2}
Let $S$ be a multiplicative subsemigroup of $\mN$. 
Then
 $$
(\mathscr{H}^\infty_S)^{-1}=\{f\in \mathscr{H}^\infty_S : \inf\limits_{s\in \mC_0} |f(s)|>0\}.
$$ 
\end{theorem}
\begin{proof} 
If $f\in (\mathscr{H}^\infty_S)^{-1}$, then there exists a $g\in \mathscr{H}^\infty_S$ 
such that for all $s\in \mC_0$, $f(s)g(s)\!=\!1$. In particular, $g\neq \bm{0}$, and so $\|g\|_\infty\!>\!0$. Thus 
$$
\textstyle 
\inf_{s\in \mC_0} |f(s)|=\inf_{s\in \mC_0} 
{\scaleobj{0.81}{\cfrac{\raisebox{-0.3em}{\text{$1$}}}{|g(s)|}}}
={\scaleobj{0.81}{\cfrac{\raisebox{-0.3em}{\text{$1$}}}{\sup_{s\in \mC_0} |g(s)|} }}
={\scaleobj{0.81}{\cfrac{\raisebox{-0.3em}{\text{$1$}}}{\|g\|_\infty}}}>0.
$$
Conversely, suppose $f\!=\!{\scaleobj{0.81}{\sum\limits_{n\in S}}}a_n n^{-s}\in \mathscr{H}^\infty_S$ is such that 
\begin{equation}
\label{15_8_2025_1757} 
\textstyle 
\delta:=\inf_{s\in \mC_0} |f(s)|>0.
\end{equation}
By \cite[Theorem~3.2]{BruSas2} (see also \cite[Theorem~2.6]{Bon}), $f\in ( \mathscr{H}^\infty)^{-1}$, i.e., ${\scaleobj{0.9}{\frac{1}{f}}}\in \mathscr{H}^\infty$. It remains to show that ${\scaleobj{0.9}{\frac{1}{f}}}\in \mathscr{H}^\infty_S$. 
Let $\epsilon\!>\!0$. 
Since we know $ \sigma_a(f)-\sigma_c(f)\!\le\!1$, and $\sigma_c(f)\!\le\! 0$, we get $\sigma_a(f)\!\le\! 1$.  Thus the Dirichlet series  given by 
$$
\textstyle 
f_{1+\epsilon}(s):=
{\raisebox{0.15em}{{\scaleobj{0.81}{\sum\limits_{n\in S}}}}}\;\! a_ n n^{-(1+\epsilon+s)}
=
{\raisebox{0.15em}{{\scaleobj{0.81}{\sum\limits_{n\in S}}}}}\;\!
{\scaleobj{0.81}{\cfrac{\raisebox{-0.21em}{\text{$a_ n$}}}{n^{1+\epsilon}}}} n^{-s}
$$
converges absolutely for all $s\in \mC$ with $\text{Re}\;\!s\!\ge\! 0$. In particular, if $s\in \mC_0$ and $\sigma:=\text{Re}\;\!s$, then 
$$
\begin{array}{rcl}
|f_{1+\epsilon}(s)-a_1|\le 
{\raisebox{0.15em}{{\scaleobj{0.81}{\sum\limits_{n\in S\setminus \{1\}}}}}}
 \;\!
{\scaleobj{0.81}{\cfrac{\raisebox{-0.21em}{\text{$|a_n|$}}}{n^{1+\epsilon}}}} \;\!n^{-\sigma}
\!\!\!\!&\le&\!\!\!( 
{\raisebox{0.15em}{{\scaleobj{0.81}{\sum\limits_{n\in S\setminus \{1\}}}}}}\;\!
|a_n|^2 )^{\frac{1}{2}}
( 
{\raisebox{0.15em}{{\scaleobj{0.81}{\sum\limits_{n\in S\setminus \{1\}}}}}}\;\!
{\scaleobj{0.81}{\cfrac{\raisebox{-0.21em}{\text{$n^{-2\sigma}$}}}{n^{2(1+\epsilon)}}}})^{\frac{1}{2}}
\\[0.21cm]
\!\!\!\!&\le&\!\!\! {\scaleobj{0.81}{\cfrac{\raisebox{-0.21em}{\text{$1$}}}{2^\sigma}}}\;\! \|f\|_2 
({\scaleobj{0.81}{\sum\limits_{n=1}^\infty}} \;\!
{\scaleobj{0.81}{\cfrac{\raisebox{-0.21em}{\text{$1$}}}{n^{2}}}})^{\frac{1}{2}}\stackrel{\sigma\to \infty}{\longrightarrow} 0. 
\end{array}
$$ 
If $a_1\!=\!0$, then \eqref{15_8_2025_1757} and the above imply $0\!<\!\delta\!\le \!0$, a contradiction. Thus $a_1\!\neq\! 0$. The above also shows that there exists a $\sigma_0\!>\!0$ such that for all $s\in \mC_{\sigma_0}$, 
we have with $\sigma\!:=\!\text{Re}\;\!s$ that 
\begin{equation}
\label{15_8_2025_1806}
\textstyle 
|f_{1+\epsilon}(s)-a_1|
=|
{\raisebox{0.15em}{{\scaleobj{0.81}{\sum\limits_{n\in S\setminus \{1\}}}}}} \;\!
{\scaleobj{0.81}{\cfrac{\raisebox{-0.21em}{\text{$a_n$}}}{n^{1+\epsilon}}}} \;\!n^{-s}|
\le
{\raisebox{0.15em}{ {\scaleobj{0.81}{\sum\limits_{n\in S\setminus \{1\}}}}}} \;\!
{\scaleobj{0.81}{\cfrac{\raisebox{-0.21em}{\text{$|a_n|$}}}{n^{1+\epsilon}}}} \;\!n^{-\sigma}
<{\scaleobj{0.81}{\cfrac{\raisebox{-0.21em}{\text{$|a_1|$}}}{2}}}.
\end{equation}
In the half-plane $\mC_{\sigma_0}$, 
$$
\textstyle 
\begin{array}{rcl} 
\quad \quad \quad
{\scaleobj{0.81}{\cfrac{\raisebox{-0.21em}{\text{$1$}}}{f_{1+\epsilon}(s)} }}
\!\!\!\!&=&\!\!\!\! (a_1+
{\raisebox{0.15em}{{\scaleobj{0.81}{\sum\limits_{n\in S\setminus \{1\}}}}}} 
\;\! {\scaleobj{0.81}{\cfrac{\raisebox{-0.21em}{\text{$a_n$}}}{n^{1+\epsilon}}}}\;\! n^{-s})^{-1}
\\[0.21cm]
\!\!\!\!&=&\!\!\!\!a_1^{-1} \pmb{(}1+{\scaleobj{0.81}{\sum\limits_{m=1}^\infty}}(-1)^m 
(a_1^{-1}
{\raisebox{0.15em}{{\scaleobj{0.81}{\sum\limits_{n\in S\setminus \{1\}}}}}} 
\;\!{\scaleobj{0.81}{\cfrac{\raisebox{-0.21em}{\text{$a_n$}}}{n^{1+\epsilon}}}}\;\!  n^{-s})^m\pmb{)},
\quad \quad \quad \phantom{Ai}(\star)
\end{array}
$$
where the geometric series converges on account of \eqref{15_8_2025_1806}. 
Since we have 
$$
\textstyle 
{\scaleobj{0.81}{\sum\limits_{m=1}^\infty}} \;\!
(|a_1|^{-1}
{\raisebox{0.15em}{{\scaleobj{0.81}{\sum\limits_{n\in S\setminus \{1\}}}}}} 
\;\! {\scaleobj{0.81}{\cfrac{\raisebox{-0.21em}{\text{$|a_n|$}}}{n^{1+\epsilon}}}}\;\!  n^{-\sigma})^m
\le 
\textstyle {\scaleobj{0.81}{\sum\limits_{m=1}^\infty}} \;\!{\scaleobj{0.81}{\cfrac{\raisebox{-0.21em}{\text{$1$}}}{2^m}}} 
<\infty,
$$
it follows that we can rearrange the terms in ($\star$), and obtain a sequence $(c_n)_{n\in \mN}$ in $\mC$ such that for all $s\in \mC_{\sigma_0}$, we have 
$$
\textstyle 
{\scaleobj{0.81}{\cfrac{\raisebox{-0.21em}{\text{$1$}}}{f_{1+\epsilon}(s)}}}
=
{\raisebox{0.15em}{{\scaleobj{0.81}{\sum\limits_{n\in S}}}}} \;\!c_n n^{-s}.
$$
Note that we used the semigroup property of $S$ here, since for $n\in S$, $(n^{-s})^m\!=\!(n^m)^{-s}$, and $n^{m}\in S$.

\noindent But if the Dirichlet series for $\frac{1}{f}\in \mathscr{H}^\infty$ is given by 
$$
\textstyle
{\scaleobj{0.81}{\cfrac{\raisebox{-0.21em}{\text{$1$}}}{f(s)}}}
=
{\raisebox{0.15em}{{\scaleobj{0.81}{\sum\limits_{n=1}^\infty}}}}\;\! b_n n^{-s}, \quad s\in \mC_0,
$$
then we obtain from the above that for $s\in \mC_{\sigma_0}$, 
$$
\textstyle 
{\raisebox{0.15em}{{\scaleobj{0.81}{\sum\limits_{n=1}^\infty}}}} 
\;\! 
{\scaleobj{0.81}{\cfrac{\raisebox{-0.21em}{\text{$b_n$}}}{n^{1+\epsilon}} }}\;\!n^{-s}
=
{\scaleobj{0.81}{\cfrac{\raisebox{-0.21em}{\text{$1$}}}{f(1+\epsilon+s)}}}
=
{\scaleobj{0.81}{\cfrac{\raisebox{-0.21em}{\text{$1$}}}{f_{1+\epsilon}(s)}}}
=
{\raisebox{0.15em}{{\scaleobj{0.81}{\sum\limits_{n\in S}}}}}\;\! c_n n^{-s}.
$$
In particular, for $n\in \mN\setminus S$, by the uniqueness of Dirichlet series coefficients (see, e.g., \cite[Theorem~7, \S5, Chap. X]{Cha}),  
$$
\textstyle
{\scaleobj{0.81}{\cfrac{\raisebox{-0.21em}{\text{$b_n$}}}{n^{1+\epsilon}} }}=0,
$$
and so $b_n\!=\!0$. This  shows that 
$$
\textstyle 
{\scaleobj{0.81}{\cfrac{\raisebox{-0.21em}{\text{$1$}}}{f(s)}}}
={\scaleobj{0.81}{\sum\limits_{n=1}^\infty}} \;\!b_n n^{-s}
={\scaleobj{0.81}{\sum\limits_{n\in S}}}\;\!b_n n^{-s},
$$
and so $\frac{1}{f}\in \mathscr{H}^\infty_S$, as wanted. 
\end{proof}

\noindent Let $\mathscr{A}_u$ be the subset of $\mathscr{H}^\infty$ 
of Dirichlet series that are uniformly continuous in $\mC_0$.
Another description of $\mathscr{A}_u$ is that 
it is the closure of Dirichlet polynomials in the $\|\cdot\|_\infty$-norm, 
see, e.g., \cite[Theorem~2.3]{ABG}. 
For a multiplicative subsemigroup $S$ of $\mN$ containing $1$, we introduce 
$
\mathscr{A}_{u,S}=\mathscr{A}_u\cap \mathscr{H}^\infty_S$. Then 
$\mathscr{A}_{u,S}$ is a unital Banach algebra with pointwise operations and the supremum norm. 

Let $\mathscr{W}$ denote the set of all Dirichlet series $f={\scaleobj{0.81}{\sum\limits_{n=1}^\infty}}  a_nn^{-s}$ such that 
$$
\textstyle \|f\|_{1}:={\scaleobj{0.81}{\sum\limits_{n=1}^\infty}} |a_n|<\infty.
$$ 
With pointwise operations and the $\|\cdot\|_1$ norm, $\mathscr{W}$ is a Banach algebra. Then 
 $
\mathscr{W} \!\subset\! \mathscr{A}_u \!\subset \!\mathscr{H}^\infty.
$ 
In the case of $\mathscr{W}$, an analogue of the classical Wiener $1/f$ lemma (\cite[p.91]{Wie}) for the unit circle holds, 
i.e., if $f\in \mathscr{W}$ is such that $\inf_{s\in \mC_+} |f(s)|>0$, then ${\scaleobj{0.9}{\frac{1}{f}}}\in \mathscr{W}$ 
(see, e.g., \cite[Theorem~1]{HewWil}, and also \cite{Edw}
% \cite{GooNew} 
for an elementary proof). 
For a multiplicative subsemigroup $S$ of $\mN$ containing $1$, we introduce 
$
\mathscr{W}_{S}=\mathscr{W}\!\cap \!\mathscr{H}^\infty_S$. Then 
$\mathscr{W}_{S}$ is a 
 unital Banach algebra with pointwise operations and the $\|\cdot\|_1$ norm. 
 We have $\mathscr{W}_S \subset \mathscr{A}_{u,S} \subset \mathscr{H}^\infty_S$.

\begin{corollary}
Let $S$ be a multiplicative subsemigroup of $\mN$ containing $1,$  and 
 $A\in \{\mathscr{A}_{u,S}, \mathscr{W}_S\}$. 
Then $f\in A^{-1}$ if and only if 
\begin{equation}
\label{18_8_2025_1842}
\textstyle 
\inf\limits_{s\in \mC_0} |f(s)|>0. 
\end{equation}
 \end{corollary}
 \begin{proof} If $f\in A^{-1}$, then $f\in (\mathscr{H}^{\infty}_S)^{-1}$. So \eqref{18_8_2025_1842} holds by Theorem~\ref{te2.2}. 
 
 Conversely, suppose $f\in A$ satisfies \eqref{18_8_2025_1842}. Then Theorem~\ref{te2.2} implies that ${\scaleobj{0.9}{\frac{1}{f}}}\in \mathscr{H}^{\infty}_S$. If $A=\mathscr{W}_S$, then by the Wiener $1/f$ theorem for $\mathscr{W}$,  ${\scaleobj{0.9}{\frac{1}{f}}}\in \mathscr{W}$, and so ${\scaleobj{0.9}{\frac{1}{f}}}\in \mathscr{W}\cap\mathscr{H}^{\infty}_S=\mathscr{W}_S$. If $A=\mathscr{A}_{u,S}$, then  \cite[Theorem~3.2]{BruSas2}  implies that ${\scaleobj{0.9}{\frac{1}{f}}}\in \mathscr{A}_u$ as well, and so ${\scaleobj{0.9}{\frac{1}{f}}}\in \mathscr{A}_u\cap\mathscr{H}^{\infty}_S=\mathscr{A}_{u,S}$. 
 \end{proof}

\section{The image of $\mathscr{H}^\infty_S$ under the Bohr transform}
 \label{section_5}

 \noindent In this section, we relate $\mathscr{H}^\infty_S$ to a natural Banach subalgebra of the Hardy algebra $H^\infty(B_{c_0})$ on the unit ball $B_{c_0}$ of $c_0$ (space of complex sequences converging to $0$ with termwise operations and the supremum norm), with vanishing derivatives of certain orders at $\bm{0}\in c_0$. 
We first introduce the following notation. 

\vspace{0.1cm}

\noindent ${\scaleobj{0.81}{\bullet}}$ $\ell^\infty$ is the Banach space of all bounded complex sequences with the 
\phantom{${\scaleobj{0.81}{\bullet}}$ }supremum norm, $\|(a_n)_{n\in \mN}\|_\infty:=\sup_{n\in \mN}|a_n|$ for $(a_n)_{n\in \mN}\in \ell^\infty$.

\vspace{0.1cm}

\noindent ${\scaleobj{0.81}{\bullet}}$ $c_0$  is the Banach space of
complex sequences tending to $0$ at infinity, 
\phantom{${\scaleobj{0.81}{\bullet}}$ }with the induced norm
from $\ell^\infty$. Let $B_{c_0}$ be the open unit ball of $c_0$.

\vspace{0.1cm}

\noindent ${\scaleobj{0.81}{\bullet}}$ $c_{00}$ is the subset of $c_0$  of all sequences in $\ell^\infty$ with compact support. 

\vspace{0.1cm}

\noindent ${\scaleobj{0.81}{\bullet}}$  Let $\bm{N}$ be the set of all compactly supported sequences that take 
\phantom{${\scaleobj{0.81}{\bullet}}$ }values in the set of nonnegative integers, i.e., $\bm{N}$ is the subset 
of $c_{00}$ 
\phantom{${\scaleobj{0.81}{\bullet}}$ }consisting of sequences whose terms belong to $\mN\!\cup\!\{0\}$. 

\noindent \phantom{${\scaleobj{0.81}{\bullet}}$ }If $\bm{\nu}\!=\!(n_k)_{k\in \mN}\in \bm{N}$ and $K\in \mN$ is such that for all $k\!>\!K$,  
$n_k\!=\!0$, \phantom{${\scaleobj{0.81}{\bullet}}$ }then we define 
$$
\begin{array}{l}
\bm{z}^{\bm{\nu}}=z_1^{n_1}\cdots z_{K}^{n_K}$ for all $\bm{z}\in B_{c_0},
\\[0.1cm]
\bm{\partial}^{\bm{\nu}}=\partial_{z_1}^{n_1}\cdots \partial_{z_{K}}^{n_K},
\\[0.1cm]
|\bm{\nu}|=n_1+\cdots+n_K,
\\[0.1cm]
\bm{\nu}!=n_1!\cdots n_K!.
\end{array}
$$

\noindent ${\scaleobj{0.81}{\bullet}}$ If $\bm{\alpha}\!=\!(\alpha_k)_{k\in \mN}, \bm{\beta}\!=\!(\beta_k)_{k\in \mN}\in \bm{N}$, then 
   $\bm{\beta}\!\preccurlyeq\! \bm{\alpha}$ if for all $k\in \mN$, $\beta_k\!\le\! \alpha_k$.

   \noindent \phantom{${\scaleobj{0.81}{\bullet}}$ }If $\bm{\alpha}, \bm{\beta}\in \bm{N}$ satisfy  
   $\bm{\beta}\!\preccurlyeq\! \bm{\alpha}$ then 
   $\textstyle 
   {\scaleobj{1.3}{{\bm{\alpha}\choose \bm{\beta}}}}
   :=
   {\scaleobj{0.9}{\cfrac{\raisebox{-0.24em}{\text{$\bm{\alpha}!$}}}{
   \raisebox{0.24em}{\text{$\bm{\beta}! (\bm{\alpha}-\bm{\beta})!$}}}}}.
   $ 

\vspace{0.21cm}

\noindent 
By the fundamental theorem of arithmetic, every $n\!\in\! \mN$ 
may be written uniquely in the form 
$$
\textstyle 
n\!=\!{\scaleobj{0.81}{\prod\limits_{k=1}^\infty}}p_k^{\nu_{k}(n)},
$$ 
where $\nu_{k}(n)\in \mN\!\cup\!\{0\}$ and $p_1\!<\!p_2\!<\!p_3\!<\!\cdots$ are all the prime numbers arranged in ascending 
 order.  We have $\bm{\nu}(n):=(\nu_k(p))_{k\in \mN}\in \bm{N}$. 

A seminal observation made by Harald Bohr \cite{Boh}, is that by putting 
$
\textstyle
z_1\!=\! 2^{-s}, 
z_2\!=\! 3^{-s}, 
z_3\!=\!5^{-s}, \cdots, 
z_n\!=\! p_n^{-s}, \cdots, 
$ 
  a Dirichlet series in
 $\sH$ can be formally considered as a power series of
infinitely many variables. Thus, the element 
 $$
\textstyle 
f(s)={\scaleobj{0.81}{\sum\limits_{n=1}^\infty}}
 a_n n^{-s} \in \mathscr{H}^\infty
$$  
gives the formal power series
$$
\textstyle 
F(\bm{z})={\scaleobj{0.81}{\sum\limits_{n=1}^\infty}}a_n 
{\scaleobj{0.81}{\prod\limits_{k=1}^\infty }}z_k^{\nu_{k}(n)},
$$
where 
$\bm{z}=(z_1,z_2,z_3,\cdots)$. We recall the precise  result below. 

 Let $H^\infty(B_{c_0})$ be the complex Banach algebra of bounded holomorphic (i.e., complex Fr\'echet differentiable) functions $F:B_{c_0}\to \mC$, with pointwise defined operations, and the supremum norm.

The Taylor series viewpoint will be useful below, and homogeneous polynomials will 
be relevant in this context. 
A function $P\!:\!c_0\!\to\mC$ is an {\em $m$-homogeneous polynomial} 
if there exists a continuous $m$-linear form $A\!:\! c_0^m\!\to\mC$, 
such that $P(\bbz) \!=\! A(\bbz,\dots,\bbz)$ for every $\bbz\in c_0$. 
The $0$-homogeneous polynomials are constant functions. 

We first recall that for a holomorphic $F:\mD^N\to \mC$, we have 
$$
\textstyle 
F(\bm{z})={\scaleobj{0.81}{\sum\limits_{m=0}^\infty}}\;
{\scaleobj{0.81}{ \sum\limits_{\bm{\alpha}\in (\mN\cup\{0\})^N,\;\!
|\bm{\alpha}|=m}}}
c_{\bm{\alpha}}(F)\bm{z}^{\bm{\alpha}}\;
\text{ for all }\bm{z}\in \mD^N,
$$
where for each $\bm{\alpha}\!=\!(\alpha_1,\cdots, \alpha_N)\in (\mN\cup\{0\})^N$, $|\bm{\alpha}|\!=\!\alpha_1\!+\cdots+\alpha_N$, and 
$$
\textstyle 
c_{\bm{\alpha}}(F)
=
{\scaleobj{0.81}{\cfrac{\raisebox{-0.3em}{\text{$1$}}}{(2\pi i)^N}}}\;\!
{\raisebox{0.12em}{$\int\limits_{{\scaleobj{0.81}{|\zeta_1|=r_1}}}$}}\cdots 
{\raisebox{0.12em}{$\int\limits_{{\scaleobj{0.81}{|\zeta_N|=r_N}}}$}}\;\!
{\scaleobj{0.81}{\cfrac{\raisebox{-0.21em}{\text{$f(\zeta_1,\cdots, \zeta_N)$}}}{
\raisebox{0.12em}{\text{$\zeta_1^{\alpha_1+1} \cdots \zeta_N^{\alpha_N+1}$} }}}}
\;d\zeta_N\cdots d\zeta_1,
$$
and arbitrary $r_1,\cdots, r_N\in (0,1)$. 
Also, 
$$
\textstyle 
c_{\bm{\alpha}}(F)=
{\scaleobj{0.81}{\cfrac{(\bm{\partial}^{\bm{\alpha}} F)(\bm{0})}{\bm{\alpha}!}}}.
$$
Then for every $m$, the function $P_m:\mC^N\to \mC$ given by 
$$
\textstyle
P_m(\bm{z})={\scaleobj{0.81}{ \sum\limits_{
\bm{\alpha}\in (\mN\cup\{0\})^N,\;\!
|\bm{\alpha}|=m}}}
c_{\bm{\alpha}}(F)\bm{z}^{\bm{\alpha}},
$$
is an $m$-homogeneous polynomial, and we have 
$$
\textstyle 
F={\scaleobj{0.81}{\sum\limits_{m=0}^\infty}} P_m\text{ pointwise on }\mD^N.
$$
The following result was shown in 
 \cite[Proposition~2.28]{DGMS}.

 \begin{proposition}
 \label{Proposition2.28DGMS}
 Let $F:B_{c_0}\to \mC$ be a bounded function. 
 
 \noindent Then the following are equivalent$\;\!:$
 
 \vspace{0.06cm}
 
 \noindent {\em (1)} $F\in H^\infty(B_{c_0})$. 
 
 \vspace{0.06cm}
 
 \noindent {\em (2)} There exists a unique sequence $(P_m)_{m\in \mN_0}$ 
of $m$-homogeneous 

\noindent \phantom{(2) }polynomials on $c_0,$ such that
$$
\textstyle 
F={\scaleobj{0.81}{\sum\limits_{m=0}^\infty}} P_m\;\text{ pointwise on }B_{c_0}.
$$
Moreover$,$ in this case$,$ we have$\;\!:$
$$
\textstyle 
{\scaleobj{0.96}{
\begin{array}{l}
P_m(\bm{z})={\scaleobj{0.81}{ \sum\limits_{
\bm{\alpha}\in (\mN\cup\{0\})^N  ,\;\!
|\bm{\alpha}|=m}}}
c_{\bm{\alpha}}(F)\bm{z}^{\bm{\alpha}} \text{ for all } \bm{z}\in B_{c_{00}},\\[-0.06cm] 
f={\scaleobj{0.81}{\sum\limits_{m=0}^\infty}} P_m\; 
\textrm{ uniformly on } rB_{c_0}\textrm{ for every } 0\!<\!r\!<\!1,\text{ and}\\[0.3cm]
\|P_m\|_\infty \!\le\! \|F\|_\infty.
\end{array}}}
$$
 \end{proposition}

The following result is known (see \cite{HedLinSei}).

\begin{proposition}
\label{prop_HLS}
The map sending elements $F\in H^\infty(B_{c_0})$ to   
$$
\textstyle 
{\scaleobj{0.96}{
f={\scaleobj{0.81}{\sum\limits_{n=1}^\infty}} {\scaleobj{0.81}{\cfrac{1}{(\bm{\nu}(n))!}}} (\bm{\partial}^{\bm{\nu}(n)}  F)(\bm{0}) 
 n^{-s},}}
$$ 
is a Banach algebra isometric isomorphism from $H^\infty(B_{c_0})$ to $\mathscr{H}^\infty$.
\end{proposition}

%Recall that for  $n={\scaleobj{0.81}{\prod\limits_{k=1}^\infty}} p_k^{\nu_k(n)} \in \mN,$ 
%$\bm{\nu}(n):=(\nu_k(n))_{k\in \mN}$. 

%\noindent 
The set $\bm{N}$ is an additive semigroup with addition given by the termwise addition of the sequences belonging to $\bm{N}$. If $\bm{S}$ is an additive subsemigroup of $\bm{N}$ containing the zero sequence $\bm{0}\!:=\!(0)_{k\in \mN}$, then let $H^\infty_{\bm{S}}(B_{c_0})$ be the subalgebra of $H^\infty(B_{c_0})$ consisting of all 
elements $F\in H^\infty(B_{c_0})$ such that 
for all $\bm{\nu}\in \bm{N} \setminus \bm{S}$, we have $(\bm{\partial}^{\bm{\nu}}  F)(\bm{0})\!=\!0$. 
The fact that $H^\infty_{\bm{S}}(B_{c_0})$ is an algebra follows immediately from the multivariable Leibniz rule, as follows.  
 If $\bm{\alpha}\in \bm{N} \setminus \bm{S}$ 
and $\bm{\beta}\preccurlyeq  \bm{\alpha}$, then either $\bm{\beta}\not\in\bm{S}$ or 
$\bm{\alpha}-\bm{\beta}\not\in \bm{S}$ (otherwise $\bm{\alpha}=\bm{\beta}+(\bm{\alpha}-\bm{\beta})\in \bm{S}$, a contradiction), and so either $(\bm{\partial}^{\bm{\beta}}F)(\bm{0})\!=\!0$ or 
$(\bm{\partial}^{\bm{\alpha}-\bm{\beta}}G)(\bm{0})\!=\!0$, showing that 
$$
\textstyle
{\scaleobj{0.96}{
(\bm{\partial}^{\bm{\alpha}}(FG))(\bm{0})
    =
   {\scaleobj{0.81}{\sum\limits_{\bm{\beta}\preccurlyeq \bm{\alpha}} }}
   {\scaleobj{1.3}{{\bm{\alpha}\choose \bm{\beta}}}}
    (\bm{\partial}^{\bm{\beta}}F)(\bm{0}) \cdot (\bm{\partial}^{\bm{\alpha}-\bm{\beta}}G)(\bm{0})
    =0,}}
$$
since each summand on the right-hand side is zero. The completeness is a consequence of the Taylor series expansion from Proposition~\ref{Proposition2.28DGMS} above. Thus $H^\infty_{\bm{S}}(B_{c_0})$ is a unital Banach subalgebra of $H^\infty(B_{c_0})$ with the supremum norm.  
% We remark that subalgebras obtained by considering vanishing derivative conditions 
% at certain points have been also considered from an algebraic perspective in the ring of polynomials, 
% see, e.g., \cite{GLTU}.

If $S$ is a multiplicative subsemigroup of $\mN$ containing $1$, then the map 
$$
\textstyle 
S\owns n \mapsto \bm{\nu}(n):=(\nu_k(n))_{k\in \mN}\in \bm{N}
$$
is an injective semigroup homomorphism, and we denote its image by $\bm{\nu}(S)$. 
An immediate corollary of  Proposition~\ref{prop_HLS} is the following. 

\begin{corollary}
\label{corollary_26_Aug}
Let $S$ be a multiplicative subsemigroup of $\mN$ containing $1,$ and let $\bm{\nu}(S)$ 
be the image of $S$ under the map $ S\owns n \mapsto \bm{\nu}(n)$. 

\noindent 
The map sending elements $F\in H^\infty_{\bm{\nu}(S)}(B_{c_0})$ to 
$$
\textstyle 
{\scaleobj{0.96}{
f={\scaleobj{0.81}{\sum\limits_{n=1}^\infty}}\;\! {\scaleobj{0.81}{\cfrac{1}{(\bm{\nu}(n))!}}}\;\! (\bm{\partial}^{\bm{\nu}(n)}  F)(\bm{0}) n^{-s},}}
$$
is a Banach algebra isometric isomorphism from $H^\infty_{\bm{\nu}(S)}(B_{c_0})$ to $\mathscr{H}^\infty_S$.
\end{corollary}

\noindent Let $A$ be a commutative unital complex semisimple Banach algebra. The dual space $A^*$ of $A$ consists of all continuous linear complex-valued maps defined on $A$. The {\em maximal ideal
space} $\text{M}(A)$ of $A$  is the set of all nonzero multiplicative elements in $A^*$ (the kernels of which are then in one-to-one correspondence with the maximal ideals of $A$). As $\text{M}(A)$ is a subset of $A^*$, it inherits the  weak-$\ast$ topology of $A^*$, called the {\em Gelfand topology} on $\text{M}(A)$. 
 The topological space $\text{M}(A)$ is a compact Hausdorff space, and is contained in the unit sphere of  the Banach space $A^*$  with the operator norm, $\|\varphi\|=\sup_{a\in A, \;\|a\|\le 1} |\varphi(a)|$ for all $\varphi \in A^*$.  
 Let $C(\text{M}(A))$ denote the Banach algebra  of 
complex-valued continuous functions on $\text{M}(A)$ with pointwise operations and the supremum norm, $\|f\|_\infty=\sup_{\varphi \in \text{M(A)}} |f(\varphi)|$ for all $f\in C(\text{M}(A))$. 
 The {\em Gelfand transform} $\widehat{a}\in C(\text{M}(A))$ of  $a\in A$ is  
 defined by $\widehat{a}(\varphi)=\varphi(a)$ for all $\varphi \in \text{M}(A)$. 

Given any $\bm{z}_*\in B_{c_0}$, the map $\varphi_{\bm{z}_*}:H^\infty_{\bm{\nu}(S)}(B_{c_0})\to \mC$ given by 
$\varphi_{\bm{z}_*}(f)=f(\bm{z}_*)$ for all $f\in H^\infty_{\bm{\nu}(S)}(B_{c_0})$ is an element of $\text{M}(H^\infty_{\bm{\nu}(S)}(B_{c_0}))$. We will use this observation in the proof of Theorem~\ref{thm_bsr} in the next and final section.

\section{Bass stable rank}
\label{section_6}

 \noindent In algebraic $K$-theory, the notion of `stable rank' of a ring
was introduced to facilitate $K$-theoretic computations 
(see \cite{Bas}). We recall the pertinent definitions below. 

Let $A$ be a unital commutative ring with unit element denoted by $1$. 
An element $(a_1,\cdots, a_n)\!\in\! A^n$ is {\em unimodular} if
there exist  $b_1,\cdots, b_n\in A$ such that
  $
 b_1 a_1+\cdots +b_n a_n=1. 
 $
The set of all unimodular elements of $A^n$ is denoted by $U_n (A)$. 
 We call $(a_1,\cdots, a_{n+1})\in U_{n+1}(A)
$ 
 {\em reducible} if there exist  $x_1,\cdots, x_n\!\in \!A $ such
that
 $
(a_1\!+\!x_1 a_{n+1},\;\!\cdots, \;\! a_n\!+\! x_n a_{n+1})\!\in\! U_n(A).
$ 
The {\em Bass stable rank} of $A$ is the least  $n\in \mN$ for
which every element in $ U_{n+1}(A)$ is reducible. The {\em Bass stable rank of $A$ is infinite} if there is no such
 $n$.
 
What is the Bass stable rank of $\mathscr{H}^\infty_S$? 

\vspace{0.03cm}

\noindent ${\scaleobj{0.81}{\bullet}}$ If $S\!=\!\{1\}$, then $\mathscr{H}^\infty_S$ is $\mC$ as a ring, and the Bass stable rank is $1$. 

\vspace{0.03cm}

\noindent ${\scaleobj{0.81}{\bullet}}$ If $S\!=\!\{p^k: k\in \mN\cup\{0\}\}$, where $p$ is a prime,  then $\mathscr{H}^\infty_S$ is isomorphic 
\phantom{${\scaleobj{0.81}{\bullet}}$ }as a Banach algebra to the classical Hardy algebra $H^\infty$,  whose Bass \phantom{${\scaleobj{0.81}{\bullet}}$ }stable rank is  $1$ (proved in \cite{Tre}).

\vspace{0.03cm}

\noindent ${\scaleobj{0.81}{\bullet}}$ If $S=\mN$, 
then the Bass stable rank of $\mathscr{H}^\infty_S=\mathscr{H}^\infty$ is infinite (shown 
\phantom{${\scaleobj{0.81}{\bullet}}$ }in \cite[Theorem~1.6]{MorSas}).

It is natural to expect that the Bass stable rank of $\mathscr{H}^\infty_S$ ought to be related to an appropriate notion of `rank/dimension' of the semigroup $S$, which perhaps gives lower or upper bounds on the Bass stable rank.
 There are several notions of the rank of a semigroup. For instance, we recall below the notion of `lower rank' and the notion of `upper rank' introduced in \cite{HowRib}. 
 For every subset $\calS$ of a semigroup $\Sigma$, there is at least one subsemigroup of $\Sigma$ containing $\calS$, namely $\Sigma$ itself. So the intersection of all the subsemigroups of $\Sigma$ containing $\calS$ is a subsemigroup of $\Sigma$ containing $\calS$, and we denote it by $\langle \calS\rangle$. For a nonempty subset $\calS$ of $\Sigma$, the subsemigroup $\langle \calS\rangle$ consists of all elements of $\Sigma$ that can be expressed as finite products of elements of $\calS$. Let $|\calS|$ denote the cardinal number of $\calS$.
 The {\em lower rank} of $\Sigma$ is defined by 
 $$
 \textstyle 
 r(\Sigma):=\inf \{|\calS|: \calS\subset \Sigma,\text{ and }\langle S\rangle=\Sigma\}.
 $$
 A subset $\calS$ of a semigroup $\Sigma$ is {\em independent} if for all $s\in \calS$, we have that $s\not\in \langle \calS\setminus\{s\}\rangle$. The {\em upper rank of $\Sigma$} is 
  $$
 \textstyle R(\Sigma):= \sup\{|\calS|: \calS\subset \Sigma, \text{ and } \calS\text{ is independent}\}.
 $$ 
 It was shown in \cite{HowRib} that $r(S)\!\le\! R(S)$.   We have the following:
 
\begin{conjecture}
Let $S$ be a multiplicative subsemigroup of $\mN$ such that $1\in S,$ and $R(S)\!=\!\infty$.
 Then the Bass stable rank of $\mathscr{H}^\infty_S$ is infinite. 
\end{conjecture}

\noindent  Let $S$ be a multiplicative subsemigroup of $\mN$ containing $1$. If $Q\subset S$ generates $S$ and $Q$ is infinite, then a cardinality argument shows that $Q$ must be countable, and so we can arrange its members in strictly increasing order as $q_1\!<q_2\!<q_3\!<\cdots$. We have the following. 

\begin{theorem}
\label{thm_bsr}
Let $S$ be a multiplicative subsemigroup of $\mN$ containing $1,$ and 
let $q_1\!<q_2\!<q_3\!<\cdots$ be a sequence in $S$ such that for all $n\in S$, 
$$
\textstyle n={\scaleobj{0.81}{\prod\limits_{k=1}^\infty}} q_k^{\alpha_k(n)} 
$$
for a unique compactly supported sequence $(\alpha_k(n))_{k\in \mN}$ of nonnegative integers.
Then the Bass stable rank of $\mathscr{H}^\infty_S$ is infinite.
\end{theorem}
\begin{proof} We follow an approach similar to the one from \cite[Theorem~1.6]{MorSas}, except that the role of the primes is now replaced by $(q_k)_{k\in \mN}$. Fix $n\in \mN$. Let $f_1,\cdots, f_{n+1}\in \mathscr{H}^\infty_S$ be given by 
$$
\textstyle 
f_1\!=\!q_1^{-s},\;\;\cdots,\;\; f_n\!=\! q_n^{-s}, \;\;f_{n+1}\!=\! {\scaleobj{0.81}{\prod\limits_{j=1}^n}}\;\! 
(\bm{1}-(q_j q_{n+j})^{-s} ) .
$$
Then $(f_1,\cdots, f_{n+1})\in U_{n+1}( \mathscr{H}^\infty_S)$ because 
by expanding the product defining $f_{n+1}$, we
obtain
$$
\textstyle 
f_{n+1}=\bm{1}-q_1^{-s}\cdot  g_1-\cdots -q_n^{-s}\cdot g_n=\bm{1}-f_1 g_1-\cdots-f_n g_n,
$$
for suitably defined $g_1,\cdots, g_n\in  \mathscr{H}^\infty_S$, and so with $g_{n+1}:=1$, we get 
$f_1g_1+\cdots+f_n g_n+f_{n+1} g_{n+1}=\bm{1}$. 
 Let
$(f_1,\cdots, f_{n+1})$ be reducible, and  $x_1,\cdots,x_n\in
\mathscr{H}^\infty_S$ be such that
$$
\textstyle (q_1^{-s}+x_1f_{n+1},\;\cdots,\;q_n^{-s}+x_nf_{n+1})\in U_n(\mathscr{H}^\infty_S).
$$ 
Let $y_1,\cdots , y_n\in \mathscr{H}^\infty_S$ be such that 
$$
\textstyle 
(q_1^{-s}+x_1f_{n+1})y_1+\cdots
+(q_n^{-s}+x_nf_{n+1})y_n=\bm{1}.
$$
Denote the isomorphism from Corollary~\ref{corollary_26_Aug} by $\iota: \mathscr{H}^\infty_S\to H^\infty_{\bm{\nu}(S)}(B_{c_0})$. Then 
$$
 \textstyle 
(\iota(q_1^{-s})+\iota(x_1)\iota(f_{n+1}))\iota(y_1)+\cdots
+(\iota(q_n^{-s})+\iota(x_n)\iota(f_{n+1}))\iota(y_n)=\bm{1}.
$$
Taking the Gelfand transform,   we obtain
$$
\textstyle
\;
(\widehat{\iota(q_1^{-s})}+ \widehat{\iota(x_1)}\widehat{\iota(f_{n+1})} )\widehat{\iota(y_1)}+\cdots+ 
(\widehat{\iota(q_n^{-s})}+\widehat{\iota(x_n)}\widehat{\iota(f_{n+1})})\widehat{\iota(y_n)}=1 .
\quad\;\;\! (\star)
$$ 
For
${\bm{z}}\!=\!(z_1,\cdots, z_n)\in \mC^n$, let $\bm{z}_*\!=\!(z_1,\cdots, z_n,\overline{z_1},\cdots, \overline{z_n},0,\cdots)\in B_{c_0}$, and 
$$
\textstyle 
\mathbf{\Phi}({\bm{z}})\!=\!\Bigg\{\!\!
\begin{array}{l} 
 -
{\scaleobj{0.81}{\prod\limits_{j=1}^n}} (1\!-\!|z_j|^2)  
{\scaleobj{0.96}{
\pmb{(}\widehat{\iota(x_1)}(\varphi_{\bm{z}_*}),\cdots\!, \widehat{\iota(x_n)}(\varphi_{\bm{z}_*})\pmb{)}
\textrm{ if } |z_j|\!<\! 1, \;\!j\!=\!1,\cdots\!,n,}}\\[0.3cm]
{\bm{0}}\; (\in \mC^n)\; \textrm{ otherwise}.
\end{array}
$$
Then $\mathbf{\Phi}$ is a continuous map from $\mC^n$ into $\mC^n$. We have that  $\mathbf{\Phi}$ vanishes outside $ \mD^n$, and so
$$
\textstyle 
 \max\limits_{{\bm{z}}\;\! \in\;\! \mD^n }\|\mathbf{\Phi}({\bm{z}})\|_2
= \sup\limits_{{\bm{z}}\;\! \in \;\! \mC^n} \|\mathbf{\Phi}({\bm{z}})\|_2,
$$
where $\|\cdot\|_2$ denotes the usual Euclidean norm in $\mC^n$. 
This implies that there must exist an $r\geq 1$ such that $\mathbf{\Phi}$
maps $K:=r \overline{\mD}^n$ into $K$. Since the set $K$ is compact and convex,
by Brouwer's Fixed Point Theorem (see, e.g., \cite[Theorem~5.28]{Rud}), it follows that there exists a
${\bm{\zeta}}\in K$ such that
 $
\mathbf{\Phi}({\bm{\zeta}})={\bm{\zeta}}.
$ 
Since $\mathbf{\Phi}$ is zero outside $\mD^n$,  we see that
${\bm{\zeta}}\in \mD^n$.  Let ${\bm{\zeta}}=(\zeta_1,\cdots,
\zeta_n)$, and $\bm{\zeta}_*=(\zeta_1,\cdots, \zeta_n,\overline{\zeta_1},\cdots, \overline{\zeta_n},0,\cdots)\in B_{c_0}$. Then for each $j\in \{1,\cdots, n\}$, we obtain
$$
\textstyle 
 \phantom{AAAAAAA}
\begin{array}{rcl}
 0\!\!\!\!&=&\!\!\!\!\zeta_j + {\scaleobj{0.81}{\prod\limits_{k=1}^n}} \;\!(1\!-\!|\lambda_k|^2) 
\widehat{\iota(x_j)}(\varphi_{\bm{\zeta}_*}) 
{\scaleobj{0.81}{\prod\limits_{k=1}^n}} \;\!(1\!-\!|\lambda_k|^2)
\\
\!\!\!\!&=&\!\!\!\! \zeta_j + (\widehat{\iota(x_j)} \widehat{\iota(f_{n+1})})(\varphi_{\bm{\zeta}_*}) .
 \phantom{AAAAAAAAAAAAAAa}
 (\star\star)
\end{array}
$$
But from ($\star$), we have
$$
\textstyle 
{\scaleobj{0.81}{\sum\limits_{j=1}^n}} ( \widehat{\iota(q_j^{-s})}+ \widehat{\iota(x_j)} \widehat{\iota(f_{n+1})}
) \widehat{\iota(y_j)}\big|_{\varphi_{\bm{\zeta}_*}}= 1 ,
$$
which together with ($\star\star$) yields $0=1$, a contradiction.  As  $n\in \mN$
was arbitrary, it follows that the Bass stable rank of $\mathscr{H}^\infty_S$ is
infinite.
\end{proof}

\noindent 
E.g., consider    
 $
\textstyle 
S=\{1\}\cup\{n: \text{there exist } x,y\in \mN\text{ such that }n\!=\!x^2+y^2\}
$ from  Example~\ref{examples}(8). Then 
 the Bass stable rank of $\mathscr{H}^\infty_S$ is infinite, as $S$ is generated by 
 $P\cup Q$, where $P$ consists of primes $p$ that are not of the form $4k +3$ for some $k\in \mN\!\cup\{0\}$, and $Q$ is the set of elements $q=p^2$, where  $p$ is a prime of the form $4k +3$ for some $k\in \mN\!\cup\{0\}$.

\end{document}